\documentclass [a4paper,14pt ]{extarticle}
\usepackage[dvips,pdftex]{graphicx,graphics}
\usepackage{amsfonts}
\usepackage[utf8]{inputenc}
\usepackage{graphicx}
\begin{document}
 
\baselineskip=16pt
\clearpage

\newcommand{\R}{\mathbb{R}}
\newcommand{\N}{\mathbb{N}}
\newcommand{\qq}{\qquad}
\newcommand{\q}{\quad}
\newcommand{\si}{\sigma}
\newcommand{\la}{\lambda}
\newcommand{\al}{\alpha}
\newcommand{\be}{{\beta}}
\newcommand{\ve}{{\varepsilon}}
\newcommand{\vk}{{\varkappa}}
\newcommand{\W}{{\bf W_k}}
\newcommand{\U}{${\bf U_{1}}$}
\newcommand{\UU}{${\bf U_{1,k}}$}
\newcommand{\Up}{${\bf U^{\,\prime}_{1,k}}$}
\newcommand{\Us}{${\bf U^{\,*}_{1,k}}$}

\baselineskip=16pt

 \setcounter{page}{1}

 MSC-2010-class:  11A25, 11N37 (Primary); 11M26 (Secondary)

\begin{center}

{\ } {\bf  \Large On Maximal Values of  Gronwall Numbers for}

\vspace{2ex}

 {\bf  \Large    Integers  with Given Greatest Prime Factor }

\vspace{2ex}

{\ } {\bf  \Large \ and Remainder in 
Modified Mertens Formula   }

\vspace{3ex}

{\bf \large by Gennadiy Kalyabin\footnote{$^)$\q Samara, 
 Russia; \ gennadiy.kalyabin@gmail.com}$^)$}

\end{center}
\vspace{1ex} 

{\it Abstract:}\ 
The new {\bf unconditional}, i.\, e. without assuming  
{\bf RH} validity, sharp
 limit relationship is found between 
the remainder  in the {\bf modified} Mertens asymptotic formula 
for the sums of primes' reciprocals and {\bf maximal} values
 $\tilde{G}_k$  of Gronwall numbers $G(N)$ among all 
integers $N$ with given greatest prime factor $p_k, \ k\to\infty,$ 
and which are multiples of the primorial $p_1\dots p_k$.

\vspace{1ex}
 The structure is described of integers at which 
corresponding maximal  values are attained.
The proofs are based on the properties of $G(N)$ studied 
in previous author's preprints.

\vspace{1ex}


\hspace{12 ex} Ramanujan-Robin inequality 

\vspace{1ex}

{\it Bibliography}:  12 items

\vspace{2ex}

{\bf 1. \ Notations, brief history  and main results}

\vspace{2ex}

\noindent
{\bf 1.1.} As usually, let $N, j, k, m,  n$ (perhaps with indices) run the set $\N$ of all positive integers,  $\N_0:=\N \cup\{0 \}$, $p$ run the set 
  $\mathbb{P}:=\{p_1, p_2, \dots \},\ p_j<p_{j+1},$ of all primes, 
$\ve$ be an arbitrary positive number, $\delta_k$ denote  sequences,
which tend to $+0$ (perhaps different even within one and the same formula); 
$\ C_y$ stand for positive constants 
which may depend only on a parameter $y$; \ symbols $\triangleright$ 
and $\Box$ denote the proof's beginning and end; \
 $\log x$ and $\gamma$ stand (resp.) for the natural logarithm
 of a positive $x$ and the Euler-Masceroni constant;
let $\theta(x), \psi(x),$ be the first and the second Chebyshev functions:
\vspace{-1ex}
$$ \theta(x):= \sum \{\log p: \ p\le x \}; 
\q \psi(x) := \sum \{\log p: \ p^m\le x \} ,
 \eqno{(1.1)}
 $$
  $$
 T(x):=\exp(\theta(x))=\prod \{p: \ p\le x\},
T_k:=T(p_k), \theta_k:=\theta(p_k),
\eqno{(1.2)}
  $$
 and let $P^+(N)$ stand for the greatest prime factor of  $N>1.$

In 1874 F. Mertens [1] proved his famous asymptotic formula
  \vspace{-1ex}  
  $$  S(x) :=\sum_{p\le x} \log \frac{p}{p-1}
  = \log \log x + \gamma+ R(x) \hbox{ with  }
   R(x)= O\left(\frac{1}{\log x}\right). 
 \eqno{(1.3)}
  $$
 
J.-L. Nicolas [2] (1983) has considered the {\bf modified} Mertens  formula\footnote{\q For $x<3 \ \log\log\theta(x)$ cannot be defined as a real number}$^)$:
$$ S(x)=\log\log \theta(x) +\gamma +Q(x), \q x\ge3,\q 
 \eqno{(1.4)}
  $$
where the remainder $Q(x):=S(x) - \log\log\theta(x)-\gamma$
is also $O({1}/{\log x})$.

\vspace{2ex}

Further let $\si(N)$ stand for the sum of all divisors of $N\in \mathbb{N}$. 
      T. Gronwall in 1913, basing on (1.3), established the  sharp  upper       order of $\si(N)$ [3]: 
     $$ \limsup_{N\to\infty} G(N)  = e^{\gamma} = 1.781\, 072\dots; 
\hbox{   where   }    G(N):= \frac{\si(N)}{N \log \log N},
  \eqno{(1.5)}
  $$
 may (and will) be  referred to as the {\it Gronwall numbers}. 

\vspace{1ex}

Denote by ${\bf W_k}$ the set of all integers $N>1$ whose greatest prime factor  $P^+(N)=p_k$, i.e. such that 
their canonical factorization into primes is:
$$ N=p_{\bf 1}^{\al_1}\cdot p_{\ 2}^{\al_2}\cdots \, p_{\ k}^{\al_k}, \q \al_k>0,
\eqno{(1.6)}
  $$
where the number $k:=k(N)$ and the exponents $\al_j\in\N_0 $
are defined uniquely, 
and   introduce the quantities:
$$G_k:=\max\{ G(N): \, N\in{\bf W_k}\, \}, \qq g_k:=\log G_k. 
\eqno{(1.7)} 
$$

Concurrently we wiil consider the set ${\bf \tilde{W}_k}$ of all those
$N\in{\bf W}_k$ which are divided by any prime 
$p_j\le p_k$, i. e. such that all $\al_j$ in (1.6) are positive 
, and the quantities
$$\tilde{G}_k:=\max\{ G(N): \, N\in{\bf \tilde{W}_k}\, \},
 \qq \tilde{g}_k:=\log \tilde{G}_k. 
\eqno{(1.8)} 
$$ 

\vspace{1ex}

{\bf Remark 1.} S. Ramanujan  (1915, the first publication  in 1997 [4]) 
and

 G. Robin (1983) [5] have established that 
{\it the validity of the inequality:

$G_k<e^{\gamma}$  (i.e $g_k<\gamma$) for all $k>4$
{\bf is equivalent} to the Riemann Hypothesis 

{\bf (RH)} on the 
non-trivial zeros of $\zeta(s)$}.

\vspace{1ex}

More detailed history of these problems may be found in [6], [7].

\vspace{1ex}

Our main goal is to establish the next  {\bf unconditional}
 limit relationship, 

which interconnects 
$\tilde{g}_k-\gamma$ and the remainder $Q(p_k)$ in (1.4).

\vspace{2ex}

{\bf Theorem 1.} {\it  Let $\ S_k:=S(p_k),
 \ \theta_k:=\theta(p_k),\, Q_k:=Q(p_k)$,
 (cf. (1.3) (1.4)). 

Then  the following limit relationship holds true: 
\vspace{-1ex}
$$ \liminf \, (\tilde{g}_k - \gamma - Q_k) {\sqrt{p_k}\,\log p_k} =  - 2\sqrt{2}, 
\eqno{(1.9)} 
$$

 \ or, in other words: 

\vspace{1ex}

 \vspace*{\fill}
 \clearpage

\ {\bf (I)}
\  there is a (constructively defined) sequence of integers

\qq  $M_k\in{\bf \tilde{W}_k}, \ k\in \N,$ such that for any
 $\ve>0$ and all $k>K_{\ve}$:
\vspace{-1ex}
$$   
\tilde{g}_k \ge \log G(M_k) > S_k -\log\log\theta_k 
- \frac{2\sqrt{2} +\ve}{\sqrt{p_k}\,\log p_k}; 
\eqno{(1.10)} 
$$

\vspace{1ex}

 {\bf (II)}   there is an {\bf infinite} set ${\bf E}\subset \N$ 
such that for any  $\ve>0$,

\qq  all $k\in{\bf E},\,k>K_{\ve}, $
  and all $N\in{\bf \tilde{W}_k}$ 
one has: } 
\vspace{-1ex}
$$
 \log G(N) \le \tilde{g}_k < S_k -\log\log\theta_k 
- \frac{2\sqrt{2} - \ve}{\sqrt{p_k}\,\log p_k}.
\eqno{(1.11)} 
$$


This assertion was  announced at the Conference dedicated 
to the 200-th P. L. Chebyshev's  anniversary [8], held 
at the Obninsk Science Center near the village
of Okatovo (his parents' estate)  in Kaluga region, 
where he was born and buried.\footnote{$^)$\q
 Pafn\'utiy Lv\'ovich Chebysh\'ev  (4[16].5.1821 -- 26.11[8.12].1894) is  the great Russian mathematician, the founder of the Saint-Petersburg mathematical school. 
He obtained fundamental results in many branches of mathematics: number theory (Edmund Landau wrote in 1909: "The first who went  the true way in the question on prime numbers and achieved  important results, was Chebyshev"), probability, uniform approximation, orthogonal polynomials etc.
Among his students were E.I. Zolotarev, G.F. Voronoy, A.M. Lyapunov, A. A. Markov (sen.).
He was elected a member of 25 academies throughout the world.  }$^)$

 The  proof of the Part {\bf (I)} of the Theorem 1
was adduced (partially) in [9].\\
In Sect. {2} the simpler construction of $M_k$ 
ensuring (1.10) is proposed.

An improved and enhanced presentation of the properties and estimates for one-step $G$-unimprovable numbers 
$N$, studied in [10], is adduced in Section {3}.
 Basing on them, the proof of the second  part 
of Theorem 1 is given in Sect. 4.

Since ${\bf \tilde{W}_k}\subset {\bf {W}_k}$ 
it is obvious that $\tilde{g}_k\le g_k.$ for all $k.$  
The (more complicated) investigation  for the sequence 
$g_k-\gamma$ will be presented in the next author's preprint.
The author hopes that this investigation will be helpful for the proof of the Ramanujan-Robin inequality.

\vspace{2ex}

{\bf 2.  \ Proof of the first part of the Theorem 1 }

\vspace{2ex}

\noindent
{\bf 2.1.  Preliminaries.} First recall that the function 
sum of divisors is {\it multi\-plicative}, i.e.
 if two naturals
 $N_1, N_2$ are mutually prime, i.e. $(N_1, N_2)=1$, then
  $\si(N_1N_2)=\si(N_1)\si(N_2)$. Taking also into account that
\vspace{-1ex}
$$ \si(p^{\al})= 1+p+p^2+\cdots p^{\al}
=\frac{p^{\al+1}-1}{p-1},
\eqno{(2.1)}
  $$
  one comes to the classical formula  for the number $N$ 
 defined by (1.6):
\vspace{-1ex}
$$ \si(N)=\si(p_1^{\al_1})\cdot \si(p_2^{\al_2})\cdots \,\si(p_k^{\al_k})
=\prod_{j=1}^k \,\frac{p_j^{\al_j+1}-1}{p_j-1},
\eqno{(2.2)}
$$
 and hence:

 \vspace*{\fill}
\clearpage
$$  \frac{\si(N)}{N} 
= \prod_{j=1}^k \frac{p_j^{\al_j+1}-1}{p_j^{\al_j}(p_j-1)} 
= \prod_{j=1}^k \left(1-\frac{1}{p_j^{\al_j+1} }\right) 
\prod_{j=1}^k \ \frac{p_j}{\ p_j-1\,}\,,
$$
$$ \qq\q \Longrightarrow \ \log \frac{\si(N)}{N} 
=\sum_{j=1}^k \log\left(1- \frac{1}{p_j^{\al_j+1}} \right)
- \sum_{j=1}^k \log\left(1- \frac{1}{p_j} \right).
\eqno{(2.3)}
$$

On the other hand, it is obvious that
\vspace{-1ex}
$$ \log N=\sum_{j=1}^k \al_j\log p_j.
\eqno{(2.4)}
$$

Joining these formulas with (1.3) and (1.4), one obtains 
$$ \log G(N) = \sum_{j=1}^k \log\left(1- \frac{1}{p_j^{\al_j+1}} \right)
 - \log\log\left(\sum_{j=1}^k \al_j\log p_j \right)
+  S_k.
\eqno{(2.5)}
$$

We will also use the following assertion establshed  in [9], Proposition 2:

\vspace{1ex}

{\bf Proposition A. } {\it Let $\la>1$; then for all 
$x> X_{\la}:=\exp(\max(1,2/(\la-1)))$ 

one has:}
$$ Y=Y(x,\la):= \sum_{p>x} \frac{1}{p^{\,\la}}
= \frac{1+\delta(x, \la)}{(\la-1) x^{\la-1}\log x}; 
\ |\delta(x, \la)|<\frac{C}{(\la-1)\log x},
\eqno{(2.6)} 
 $$ 

 {\bf 2.2. } Next the explicit construction of $M_k$ 
is adduced,  which gives 

\ the proof for the part {\bf(I)} of the Theorem 1.
\vspace{1ex}

{\bf Lemma 1.} {\it Let us choose three sequences of naturals 
$n=n_k, \,  m=m_k,$ 

\ $b=b_k$  such that for some $\ve_0,\,0<\ve_0<0.5,$ and $k\to\infty$:
$$ {\bf(i)}\ \theta_n\approx \sqrt{2\theta_k}\,,
 \, {\bf(ii)}\ \theta_n=o(\theta_m^2),
\, {\bf(iii)}\ b_k >\frac{{\log\theta_k}}{\log2-\ve_0}\, , 
\, {\bf(iv)}\ b_k\, \theta_{m} = o(\theta_n). 
\eqno{(2.7)} 
 $$ 

 Then  the naturals 
\vspace{-2ex}
\vspace{-1ex}
$$ M_k:=\prod_{j=1}^m p_j ^{b_k}
\ \cdot \prod_{j=m+1}^n p_j^2
\ \cdot \prod_{j=n+1}^k p_j ,
\eqno{(2.8)} 
 $$ 

\vspace{-2ex}
 \ satisfy (1.10).} 

\vspace{1ex}

 $\triangleright$\ First note that by virtue of the
 relationships
$\theta_k\approx p_k \approx k\log k$, 

the conditions (2.7) are
compatible: one may take, e. g.:
 $$ n_k:=\bigl[ \sqrt{8k\,/\log k\,}\,\bigr], 
\q m_k:= \bigl[\,k^{1/3}\,\bigr], 
\q  b_k: = [\,1.5\log k\,],
\eqno{(2.9)} 
 $$

\vspace*{\fill}

\clearpage

\q\ From (2.4), (2.7) and (2.8) it follows immediately that:
$$  \log M_k = b_k\,\theta_m +2\,(\theta_n-\theta_m)+(\theta_k-\theta_n)
$$
$$=\theta_k +\theta_n +(b_k-2)\,\theta_m = \theta_k+\sqrt{\,2\theta_k\, }
+ o ( \sqrt{ \theta_k }\,)
\eqno{(2.10)}
$$

\q\ whence, with Taylor formula one obtains in turn  that 

\q\ for certain  $x\in (\theta_k, \log M_k)$:
\vspace{-2ex}
$$ \log\log \log M_k = \log\log\theta_k
+\frac{\sqrt{2}}{\sqrt{\theta_k}\log \theta_k} 
- \frac{\,2\theta_k\, (\log x +1)}{x^2\,\log^2x}.
\eqno{(2.11)}
$$

\q\ Analogously, from (2.3), (2.7), (2.8) one obtains (explanations below):
\vspace{-1ex}
$$ \log\frac{\si(M_k)}{M_k}= \sum_{j=1}^m\log\left(1-\frac{1}{p_j^{b_k+1}}\right) 
+ \sum_{j=m+1}^n\log\left(1-\frac{1}{p_j^3}\right)
$$
$$ \q\ + \sum_{j=n+1}^k\log\left(1-\frac{1}{p_j^2}\right) 
+S_k
$$
$$
> S_k - \frac{ b_k m}{ 2^{b_k+1}} 
-   \frac{1+o(1)}{2p_m^{2}\log p_m}
 -\frac{1+o(1)}{p_n\log p_n} 
= S_k - \frac{\sqrt{2}+o(1)}{\sqrt{p_k}\log p_k} .
\eqno{(2.12)}
$$

\vspace{1ex}

\qq Here the Proposition A is applied to estimate the second  and the third 

\q\ sums in (2.12) with $\la=3,\,\la=2$ (resp.), and  
the relationships are used, 

\q\ which follow from  (2.7): 
 $ \theta_n=o(\theta_m^2),\, b_km_k\sqrt{p_k} \,\log p_k=o(2^{b_k}),$

\q\ $\ p_n\log p_n \approx
 \sqrt{p_k}\, \log p_k / \sqrt{2}$.

\vspace{1ex}

\qq Joining (2.11) and (2.12) one comes to (1.10) $\Box$

\vspace{2ex} 
  
{\bf 3.  \ Diverse conditions of  $G$-maximality }

\vspace{2ex}

{\bf 3.1. }  G. Caveney, J.-L. Nicolas and J. Sondow [6], [7] 
have introduced the 

classes of numbers {\bf GA1}, {\bf GA2}.

\vspace{1ex}

{\bf Definition 1} (cf [7], Sect. 2). {\bf (i)}:
{\it  An integer $N$ belongs  to  ${\bf GA1}$ if it is

 composite and for any prime factor $q$ of $N$ one has:} $G(N/q)\le G(N)$. 

\vspace{1ex}

{\bf (ii)}: {\it  An integer $N$ belongs  to  ${\bf GA2}$ 
if $G(Na)\le G(N)$ for any integer $a$}. 

\vspace{1ex}

{\bf (iii)}: {\it  An integer $N$ is called {\bf extraodinary} 
if it is both in  ${\bf GA1}$ and  ${\bf GA2}$}.

\vspace{1ex}

Each of these classes are not empty and ${\bf GA1}$ is infinite.

The cardinality of ${\bf GA2}$ is not known and this is not accidental.

\vspace{1ex}

{\bf Proposition CNS} [6]. {\it {\bf RH} is equivalent to each of the conditions:}

 {\bf(i)}:  4 {\it is the {\bf only extraordinary}  number;  }
 {\bf(ii)}:  $\#{\bf GA2} < \infty;$
  
{\bf(iii)}:  $N>5040 \Rightarrow N\notin{\bf GA2}$.

\vspace{1ex}

\vspace*{\fill}

\clearpage

{\bf Remark 2.} The author [10] has considered the class ${\bf U_1}$ of {\it one-step 

G-unimprovable}
 numbers which differ from extraordinary ones by replacing 

a difficult to verify condition  {\bf GA2} by more  constructive  relationship 

involving the   
multiplication by {\bf single primes} only. 

The 
infinitude of  ${\bf U_1}$ was established, 
the least numbers in it
being 

 $N^*_1=2\cdot7=14, \ 
N^*_2=T(23)\cdot T(5)\cdot T(3)\cdot 2^2 $
$=160\, 626\, 866\, 400.$

\vspace{1ex}

{\bf Definition 2} (cf [10], Sect. 1). {\it  An integer $N$ belongs
 to  ${\bf U}_1$ if: 

{\bf (i)} it belongs to {\bf GA1} and
 {\bf(ii)}  \    $G(Np) \le G(N)$  for any {\bf prime} $p$.}

 \vspace{1ex}

We will also consider  the local versions of this class.

\vspace{1ex}

{\bf Definition 3} (cf [10], Sect. 2). {\it  An integer $N\in{\bf W_k}$ 
is {\bf locally $G$-maximal } 

$(N\in{\bf U_{1,k}})$ if:
\ {\bf (i)} it belongs to {\bf GA1} and
 
{\bf(ii)}  \    $G(Np) \le G(N)$  for any prime $p\le p_k$.}

\vspace{1ex}

For infinitely many $k$ there are  no integers 
$N\in {\bf U_{1,k}}$. So we introduce 

\vspace{1ex}

{\bf Definition 4.} {\it Let 
{\bf E}$:=\{k\in\N: {\bf U_{1,k} }\not=\emptyset\}$. }

\vspace{1ex}

It has been  proved in [10] that the sets {\bf E} and 
$\N\setminus {\bf E}$ are both infinite.

\vspace{1ex}

Concerning the {\bf constructiveness} of this definition and  {\bf infinitude} of {\bf E} 

cf. Prop. 5{\bf(IV), (V)} and Remark 5, 6,\ S. 3.5  below.

 \vspace{1ex}

{\bf 3.2.} In this Section  some helpful charaterizations of numbers
 in ${\bf U}_{1,k}$ and 

${\bf U}_{1}$ are established (Theorem 2), which  
will be essentially used in the proof 

of the second part of Theorem 1.

\vspace{1ex}
 We will  need some more notations and  auxillary assertions.

 \vspace{1ex}

The notation $ N \ \|\ p^{\al},\ p\in \mathbb{P},\ \al \in \N_0$
 means that    $ N=p^{\al}m, \ (m, p)=1$.

For $ p>1,\ \eta>1,\ \al\ge0$ introduce the quantities:
\vspace{-1ex}
$$   \nu = \nu(p, \eta):=\frac{1}{\log p} 
\log\left(\frac{(p-1)\log\eta}
{p^2\log\left(1+\frac{\log p}{\eta}\right)}+\frac{1}{p}\right);
\eqno{(3.1)}
$$
$$
\lambda = \lambda(p,\alpha):=
 \frac{p^{\alpha+2} -1}{p^{\alpha+2} -p}
=1+\frac{1}{p+p^2 + \dots +p^{\alpha+1}};
 \eqno{(3.2)}
$$
and let  $\xi:=\xi(p, \alpha)\  $ be the unique positive
root of the equation:
 $$ \xi^{\lambda(p, \alpha)}= \xi+\log p
\iff \log \xi =\frac{p^{\alpha+2} -p\ }{p-1}
 \ \log\left(1 + \frac{\log p}{\xi}\right) .
 \eqno{(3.3)}
 $$

From the very form of the equation  it becomes clear that its root
$\xi(p,\alpha)$ increases monotonically as  
 $\alpha,\ \alpha \ge 0,$ increases. 

\vspace{2ex}

\vspace*{\fill}

\clearpage

{\bf Proposition 1.} 
{\it Let  $N\  \| \ p^{\alpha}, \ N>2, \ \eta:=\log N$;  then the following
four  conditions are equivalent} (cf. Def. 3{\bf(ii)}): 

\vspace{-1ex}
 $$ \qq\q \hbox{\bf (i)}\ G(Np) \le G(N),\q
\hbox{\bf (ii)}  \q \eta^{\lambda} \le \eta +\log p,\q
$$
$$
 \qq  \hbox{\bf (iii)}\q \eta\le \xi(p, \alpha), 
\qq   \hbox{\bf (iv)}\q \al\le \nu(p, \eta). \q
\eqno{(3.4)}
 $$

\vspace{1ex}

 The equivalence of {\bf (i)}, {\bf (ii)} and {\bf (iii)} in (3.4),
as well as Propositions 2 and 3 below, have been proved in 
 [10], Sect. 2.

\vspace{1ex}

\noindent
{\bf Remark 3.}
It is interesting to note 
that the {\it explicit}
 (though cumbersome) inequality $\al\le\nu(p, \eta)$ 
is equivalent to the {\it implicit}
 relationship $\eta\le\xi(p,\al)$. 

\noindent
$\triangleright$ In fact, by virtue of defining
formulas (3.1), (3.3) and monotone increasing of the left-hand  
and decreasing if the right-hand side of the last equation (3.3)
(with respect to $\xi$) one obtains the following chain of equivalences:
\vspace{-1ex}
$$\al \le \nu(p,\eta) 
\iff
p^{\al} \le \frac{(p-1)\log\eta}
{p^2\log\left(1+\frac{\log p}{\eta}\right)}+\frac{1}{p}
 $$
\vspace{-1ex}
$$  \iff \frac{p^{\al+2}-p}{p-1} \ 
\log\left(1+\frac{\log p}{\eta} \right)
\le \log\eta \iff \eta \le \xi(p,\al)\ \Box
$$

\vspace{1ex}

{\bf Proposition 2.}
{\it Under assumptions and notations of previous assertion,  
let $N/p>2,\ \al>0$;  then the following four  conditions are equivalent: }
\vspace{-1ex}
 $$\q  \hbox{\bf (i)}\q G(N/p) \le G(N),\q
\q \hbox{\bf (ii)}  \q (\eta-\log p)^{\lambda} \ge \eta ,\qq\q
$$
\vspace{-2ex}
$$
\qq\qq   \hbox{\bf (iii)}\ \eta\ge \xi(p, \al-1)+ \log p,
\q   \hbox{\bf (iv)}\q  \al\ge 1+\nu(p, \eta-\log p).
\eqno{(3.5)}
 $$

\vspace{1ex}

{\bf Proposition 3.}
{\it For all  $\alpha\in \N_0, \ p>1
$ the following inequalities hold:}
\vspace{-1ex}
$$ \hbox{\bf (i) }   p  -   { \log p}  <\xi(p,0)< p ; \ 
\hbox{\bf (ii) } \   \frac{p^{\alpha+1}}{ \alpha +1}
  <\xi(p,\alpha) <  C_p\ \frac{ p^{\alpha+1}}{ \alpha +1}, 
\ \al>0;
\eqno{(3.6)}
$$

{\it   where  }
\vspace{-1ex}
$$    
\ C_p:=\left\{\qq 3; \qq \q p =2, 
\atop \frac{p \log p}{(p-1)(\log p -1/e)}; 
\ p\ge 3 \right.  
\quad \alpha\in \mathbb{N}; \q C_p\to 1, \ p\to\infty.
 \eqno{(3.7)}
$$

\vspace{1ex}

Certainly we will need the classical estimates 
for Chebyshev functions.

\vspace{1ex}

{\bf Proposition 4.}  {\it The  following relationships hold} 
(cf [12], p. 111, 131):
$$ \theta(x)\le \psi(x)<1.04 x,\, x\ge2;
\eqno{(3.8)}
$$
$$ {\bf(i)}\ \theta(x)\approx \psi(x) \approx x,\q {\bf(ii)}
\ \psi(x) - \theta(x) \approx\sqrt{x}, \q x\to +\infty;
\eqno{(3.9)}
$$
$$\ {\bf(i)}\ p_k\approx k\log k; 
\q {\bf(ii)}\ k\approx p_k/ \log p_k,   \ k\to +\infty.\qq
\eqno{(3.10)}
$$

\vspace*{\fill}
\clearpage

\noindent
{\bf 3.3.} Now we can outline the properties  
of numbers from \UU\ and \U .

\vspace{1ex}
 
{\bf Theorem 2.} {\it\ Let $k>4,$\footnote{$^)$ The values $k\le 4$ are
excluded because, e. g. the number $N^*_1:=14=2\cdot7$ 
belongs to ${\bf U}_{1}$, 

\qq  but the monotonicity condition (3.11){\bf(i)} is violated:
 $\al_2(14) =\al_3(14)=0$.  }$^)$ 
 $ N\in\W$ be  given by (1.6), 

$\eta:=\log N=\sum_{j=1}^{k} \al_j\log p_j$, $\ \nu(p, \eta),$ 
and $ \xi(p,\al)$ be  defined by (3.1), (3.3).}

\vspace{1ex}

\noindent
 {\bf(I)} 
 {\it  If $N$ belongs to  \UU \,,  and thus $k\in{\bf E}$,}
 (cf Def. 4), {\it  then one has:}
\vspace{-1ex}
$$ \q \ \ {\bf (i)} 
{ \ the \ exponents' \ monotonicity}: 
\al_1\ge \al_2  \dots \ge \al_k; 
 \ {\bf (ii)}\  \al_k=1;
\eqno{(3.11)}
$$

\vspace{-1ex}
{\bf (iii)}\ {\it  the estimates}
\vspace{-1ex}
$$ 
 \xi(p_j,\alpha_j -1) +\log p_j 
 \le \eta \le  \xi(p_j,\alpha_j), \ \forall j\le k; 
\eqno{(3.12)}
$$

{\it  hold true, or equivalently,}
\vspace{-1ex}
$$ 
 1+\nu(p_i, \eta-\log p_i)\le \al_i \le \nu(p_i, \eta); \ \forall i\le k; 
\eqno{(3.12^{\,\prime})}
$$

{\bf (iv)}\ {\it   for any $\ve>0$ there is $K_{\ve}\in\N$ such that }
\vspace{-1ex}
$$ k>K_{\ve}, \ k\in{\bf E}\ \Longrightarrow
\ \bigl| \eta - \theta_k  -  \sqrt{\, 2 \,\theta_k\,}\bigr| 
< \ve \sqrt{\theta_k}, 
\eqno{(3.13)} 
$$
\vspace{-1ex}
{\it or, in other words, there is such a  sequence 
$\{\delta_k\}_{k\in{\bf E}} \searrow 0,$ that}:
\vspace{-1ex}
$$ \eta = \theta_k + B(N) \sqrt{\theta_k}; 
\q  \sup\{\bigl| B(N) -\sqrt{2}\bigr|: N\in {\bf U_{1,k}}\}<\delta_k.
\eqno{(3.13^{\,\prime})} 
$$

\noindent
{\bf(II)} {\it Conversely, if  $N\in {\bf W_k}$ and
the relationships (3.12)  are fulfilled, then

\q $k\in{\bf E}, \ N\in$\, \UU ,  and whence
 (3.11),  (3.13) and   (3.1$3^{\, \prime}$)   also hold true.}

\vspace{1ex}

\noindent
{\bf(III)} {\it $N\in$\,\U\ iff $N\in$\,\UU\ } (cf Def. 2, 3),
{\it and $\eta\le\xi(p_{k+1}, 0)$}.

\vspace{2ex}

\noindent
$\triangleright$ 1)  First  show that 
{\it if $N\in{\bf U_{1,k}},\,  k>4,$ 
then all the exponents}  $\al_j\ge1,\ $ 
(otherwise the quantity 
$\xi(p_j, \al_j-1)$ in (3.12) would  be undefined). 

We start with obvious inequality:
 $\eta=\log N\ge \log p_k\ge \log p_5 = \log11$; 
 recalling that (cf. (3.2)) $\la(p,0) =1+1/p$, one obtains
 $\eta^{\la(2,0)} > \eta +\log2$, whence
 by virtue of Proposition 1, cf. (3.4){\bf(ii)}, it follows that  $\al_1\ge1,$
and thus
 $N=2\cdot11\cdot n_k^{(1)}\ge22;\, (n_{\,k}^{(1)}, p_i)=1\ \forall i>k$. 
Analogously,  $\eta\ge\log22=3.091..$  implies 
 $\eta^{1/3}\ge 1.456... > 1+(\log3)/\eta =1.355... $,
 whence $\eta^{4/3}>\eta+\log3$, which means (again by virtue
of (3.5){\bf(ii)}), that  $\al_2\ge1$. (This argument  doesn't
work out  for the number $N_1^*=14\in{\bf U}_1$).

Further, presenting $N $ as $N=p_k n_k$ with  
$(n_k, p_i)=1\ \forall i>k $, one has $G(n_k)\le G(N)$, because otherwise $N$ wouldn't satisfy Def. 3{\bf(i)}; hence according to Proposition 1 
$\log n_k > \xi(p_k, \al_k-1)$. 

Supposing that some exponent $\al_j=0,\ 2\le j<k,$ and taking into account the monotonicity properties of $\xi(p,\al)$, one comes to: 
$$ \log N>\log n_k \ge \xi (p_k, \al_k-1)> \xi (p_j, 0);
\eqno{(3.14)}
$$
and    applying  the Proposition 1 once again, one obtains $G(Np_j)>G(N)$,
which contradicts to Def. 2{\bf(i)}, Def. 1{\bf(i)}$\ \Box.$

\vspace*{\fill}
\clearpage

\vspace{1ex} 

\noindent
$\triangleright$  2) Taking into account  the Definition  2, 
one concludes that  inequalities (3.12) are  merely 
reformulations of Propositions 1, 2 $\Box$.

\vspace{1ex} 

\noindent
$\triangleright$  3) 
Suppose now that monotonicity condition 
(3.11){\bf (i)}  is violated; this means that $\al_j<\al_i$ 
for some $j<i\le k$; then by virtue of monotonic increase of
   $\xi(p,\al)$  with respect to  $p$ and $\al$,  we obtain 
$ \xi(p_i, \alpha_i-1)+\log p_i \ge
 \xi(p_i, \alpha_j) +\log p_i  >  \xi(p_j, \alpha_j).$
Therefore under this assumption the  system of inequalities (3.12){\bf(i)} 
for $\eta:=\log N$  would be incompatible, 
and thus $N\notin {\bf U}_{1,k} \ \Box$.

Note that here we used Proposition 2, and thus meant that $N>2p_k$.

So the above reasoning is not applicable to the number $N^*_1=14=2\cdot7$.

\vspace{1ex} 

\noindent
$\triangleright$\ 4) Now we proceed to the proof of the estimate (3.13). 

For any $N\in{\bf W_k},$ (cf. (2.1)),  
 satisfying the monotonicity condition (3.11){\bf (i)}, and any
 $ m\in \{1, 2, \dots, \alpha_1(N) \}$, let us introduce  the integers
\vspace{-1ex}
 $$ q_m=q_{m}(N):= \max\{p_j: j\le k, \alpha_j\ge m\},\, 
q_1=p_k\ge q_2\ge\dots\ge q_{\al_1}=2. 
\eqno{(3.15)}
$$  

Then  by virtue of (1.2) one has:
\vspace{-1ex} 
$$N=\prod_{m=1}^{\al_1} T(q_m),\qq
\eta=\sum_{m=1}^{\alpha_1} \theta(q_m);
 \eqno{(3.16)}
$$  

Therefore,  if $N\in{\bf U_{1,k}}$ then from the inequalities: 
$\xi(p_j,\alpha_j -1) < \eta $, 
(cf.\\ (3.6){\bf (ii)}), (3.4){\bf (iii)}) and  Proposition 3, 
it follows that $q_m < (m \eta)^{1/m},$
for any $\,m\ge2$, and taking into account (3.8) 
and the relationship
$m^{1/m}\le3^{1/3}= 1.442..\, \forall m\in\N$,
 one comes to the upper estimate 
 $\theta(q_m)< 1.516\, \eta^{1/m}$.

Hence, from (3.16)
one obtains for all $\, k\in{\bf E},\,k>k_0$ (explanations below):
\vspace{-1ex}
$$ \eta < \theta_k +\theta(q_2) 
+1.516 \eta^{1/3}  \al_1 
< \theta_k 
+ (\sqrt{2}+\delta_k)  \sqrt{\eta}+   \eta^{1/3}\log\eta
$$
$$<\theta_k+(\sqrt{2}+2\delta_k)  \sqrt{\eta}. 
\eqno{(3.17)}
$$

Here we have also taken into account that by virtue 
of right inequality (3.1$2*{\,\prime}$) with $i=1$ 
and defining formula (3.1):
\vspace{-1ex}
$$ \al_1\le\nu(2,\eta) =  \frac{1}{\log 2} 
\log\left(\frac{\log\eta}
{4\log\left(1+\frac{\log 2}{\eta}\right)}+\frac{1}{2}\right)
< A_1+A_2\log\eta,
\eqno{(3.18)}
$$  
where $A_1, A_2$ are certain positive absolute constants.

Now putting
 $\eta=\theta_k+B\sqrt{\theta_k}, B\ge0,\, A:=\sqrt{2}+2\delta_k$
and substituting it into (3.17), one comes
 (after elementary transformations) to:
\vspace{-1ex}
$$ B < A \sqrt{1+\frac{B}{\sqrt{\theta_k} } }
< A\left(1+\frac{B}{2\sqrt{\theta_k}}\right)
\Rightarrow B <\sqrt{2}+3\delta_k, 
\ \forall k>k_1. 
\eqno{(3.19)}
$$  

\noindent 
$\triangleright$\ 5) To prove that the {\it lower} limit of 
$B(N), \ k\in {\bf E},$ also equals  $\sqrt{2}$,
 let us note that for $k\to\infty$ according to 
the right inequality (3.1$2^{\,\prime}$) with $i=2$,  defining formula (3.16),
 relationships (3.7){(ii)} with $\al=1$ and (3.8),  one has
 $|\theta(q_2(N)) - \sqrt{2\theta_k}| < \delta_k\to0;$ whence it follows 
$\eta(N)-\theta_k> \sqrt{2\theta_k} (1-\delta_k)$.
 Therefore one may assert that 
$\min\{\log N: N\in{\bf U_{1,k}}\}
 - \theta_k\approx \sqrt{2\theta_k} $, 
which jointly with the upper estimate (3.17)
leads to
 $\sup_{N\in{\bf U_{1,k}}}\,|B(N) - \sqrt{2}|<\delta_k\ \Box $.

\vspace{1ex} 

\noindent
$\triangleright$ 6) \ Futher, assuming  $\al_k>1$ and applying
the monotonicity condition (3.11){\bf (i)}, 

\ one comes to 
$\eta\ge2\theta_k$, which would contradict to the estimate (3.13). 

\ Thus for any $N\in{\bf U}_{1,k}$ 
the relationship (3.11){\bf (ii)}
necessarily holds  $\ \Box$. 

\vspace{1ex} 

\noindent
$\triangleright$ 7) \ The condition (3.11) jointly with (1.2) implies $N>2p_k$, 
whereas the right and the left inequalities (3.12), by virtue of
of Propositions 1 and 2, coincide  (resp.) with the conditions
{\bf(i), (iii)} of Definition 2.

\vspace{1ex} 

\noindent
$\triangleright$ 8) \ Moving to the Part {\bf(II)}, let us
suppose that $\al_j\ge1\ \forall j\le k$; then by virtue of Propositions 
1 and 2, the system of inequalities (3.12) is equivalent
 to the relationships
\vspace{-1ex}
$$G(N)\ge \max(G(Np_j), G(N/p_j)), \, \forall \,j\le k.
\eqno{(3.20)}
$$

\vspace{-1ex}

Further,  (3.12) and the increase $\xi(p,\al)$ with respect to $p$
imply that for any $i>k$ one has:
 $\xi(p_i,0)\ge\xi(p_{k+1},0)>\eta$,
 i\,. e. $G(Np_i)<G(N), \, i>k$,
which jointly with (3.20) and Definition 2 
means that $N\in{\bf U_1}$. 

Thus the proof of the  Theorem 2  is complete $\ \Box.$ 

\vspace{1ex} 

\noindent
{\bf 3.4.}\  Denote by ${\bf U^{\,\prime}_{1,k}}, k>4 $ the set 
of all integers $N\in {\bf U^{\,\prime}_{1,k}}$ which divide all 
$p_j,\,j\le k$, i. e. satisfy (3.20). 
 By definitions and from the above reasonings (cf. left inequlity (3.12)) it follows easily that
 {\it an integer $N \in {\bf W_k},\,k\ge5,$ 
belongs to ${\bf U_{1,k}}$ iff 
 \ $N\in{\bf U^{\,\prime}_{1,k}}$ 
and} $ \xi(p_{k}, 0)+\log p_k \le \eta .$

\vspace{1ex} 

{\bf Remark 4.}
All the classes ${\bf U^{\,\prime}_{1,k}}$ in contrast to ${\bf U_{1,k}}$
are not empty.

  Next the algorithm, 
proposed in [10, th.2, 4], is described
 which yields 
 the 

{\it minimal}
element
 $V_k\in{\bf U}^{\,\prime}_{1,k}, \ k>4,$
among which infinitely many $V_k\in{\bf U_1}$.

 {\bf Proposition 5.}  {\it Let $k>4$};{ put  $ V_{\ k}^{(0)} =T(p_k) :=p_1\cdot  \dots \cdot p_k$ (cf. (2.1)), and for $s\in \mathbb{N}$   define  inductively (cf. (3.15))  the integers:
\vspace{-1ex}
$${\beta}_{j, k, s}:= \max \{\beta\in \mathbb{N}:
\ \xi(p_j,\beta -1)+\log p_j \le \log V_{\ k}^{(s-1)} \};\ j<k;
$$
\vspace{-2ex}
$$
{\beta}_{k, k, s}:= 1; \qquad  V_{\ k}^{(s)} :=  \prod_{j=1}^k \ p_{\ j}^{\beta_{j, k, s}}.
\eqno{(3.21)}
$$

\vspace{-1ex}

{\it Then}:

 {\bf (I)}\ {\it for fixed $j, k$ the numbers  
$\be_{j, k,s},  V_{\ k}^{(s)},$ 
as well as  $G(V_{\ k}^{(s)})$ do not decrease as   $s$ increases and
are bounded from above; hence there is  $s_0:=s_0(k)$ such, that for all 
 $s\ge s_0 $ the relationships  $ V_{\ k}^{(s)} =:V_k, 
\ \be_{j, k,s}=\be_{j,k} $ hold.

\vspace*{\fill}

\clearpage

For stabilized values the inequalities are fulfilled}:
\vspace{-1ex}
$$\xi(p_j,{\be}_{j, k} -1) +\log p_j \ 
 \le \log V_k <  \xi(p_j, {\be}_{j, k} ) ; 
\qquad 1\le \  j\ \le k;  
\eqno{(3.22)}
$$

{\it which by virtue of Theorem 2$^{\,\prime}$, part {\bf(II)}, means, that 
$V_k\in{\bf U^{\,\prime}_{1,k}}$.}
\vspace{1ex}

{\bf (II)}\ {\it For any other integer 
${N}:=p_{\,1}^{{\alpha}_{1}}\cdots p_{\, k}^{{\alpha}_{k}}
\in {\bf U^{\,\prime}_{1,k}}$ one has: $\al_{j}\ge\be_{j,k}$,

\hspace{5ex} $\forall j\le k$, \ i.\,e.\, $V_k\, |\, N,$
 and thus $V_k$ is the least element in \Up .}
\vspace{1ex}

{\bf (III)}\ {\it Every $V_k, \ k>4, $ is a divisor of any} $\ V_{m}, m>k $. 

\vspace{2ex}

\,{\bf (IV)}\ {\it There are infinitely many indices $k=k_m\nearrow\infty$
at which the sequence  

\ $G(V_k)$ has local maxima, \,i.\,e.
$G(V_{k_m})\ge \max(G(V_{k_m-1}, G(V_{k_m+1}))$, 
and each 

\ of  these integers $V_{k_m}\in{\bf U_1}$; 
 thus $ \{k_m\}_1^{\infty}\subset {\bf E}$, cf. Def. 2,\, 3,\, 4.} 

\vspace{1ex} 

\ {\bf (V)}\ $\log V_k \ge \xi(p_k, 0)+\log p_k
 \iff {\bf U_{1,k} }={\bf U^{\,\prime}_{1,k}}
 \not=\emptyset\ \iff k\in{\bf E}.$ 

\vspace{2ex} 

{\bf Remark 5.} According to J.E. Littlewood result 
 (1914, cf [12, p. 322])  one has: 
$\theta(x) -x = \Omega_{\pm} \sqrt{x}\ \log\log\log x $.
Therefore there are {\bf infinitely many} $k\in \N$ such that 
$\theta_k>p_k - \sqrt{p_k}$. By virtue of Th. 2 all such 
$k\in{\bf E}$.

\vspace{1ex}

{\bf Remark 6.} To summarize  the results stated 
above one may draw 
  the following relationships between the classes defined:
$$ {\bf U_1}\cap{\bf W_k} \subset {\bf U_{1,k} } 
\subset {\bf U^{\,\prime}_{1,k} }\not=\emptyset. 
\eqno{(3.23)}
$$

\vspace{2ex} 

{\bf 4. \ Proof of the second part of Theorem 1.}

\vspace{2ex} 

\noindent
{\bf 4.1.}\ We have and intend 
to prove that {\it if $k\in{\bf E},\, k>4,$ then for any $N\in{\bf W_{k}}$
one has} (cf. (1.11)): 
\vspace{-2ex} 
$$ \log G(N)< S_k - \log\log\theta_k 
- \frac{2\sqrt{2} - \delta_k}{\,\sqrt{p_k}\log p_k}, \q \delta_k\searrow0.
\eqno{(4.1)}
$$

But according to Theorem 2\,{\bf (I)},  for $k\in{\bf E}$
the maximum $G(N): N\in{\bf W_k}$ is attained at the integer
$Y_k\in {\bf U_{1,k}}$. So it will be  sufficient 
to establish (4.1) for $N\in {\bf U_{1,k}}, \,k\in{\bf E},$ {\bf only.}  

 We will continue to take use of the notations and identities 
(2.1) -- (2.5)  of the Sect.  2.1. 
Let us fix any   $N\in{\bf U_{1,k}},\ k>4;$ then by virtue 
of Theorem 2   the exponents $\al_j $ do not increase; 
denote by $n, m$ the maximal $j$ such that $\al_j\ge3$ or
 (resp.) $\al_j\ge2$, and consider a function of $m$ 
"liberated" real exponents $t_j>0, \ 1\le j\le m$:

\vspace{-2ex} 
$$  \tilde{N}(t) = \tilde{N}_{n,m}(t_1, t_2,\dots, t_m):= 
\prod_{j=1}^m p_{\ j}^{t_j}\, \cdot \prod_{j=m+1}^n p_{j}^{2}
\, \cdot \prod_{j=n+1}^k p_j \, ,
\eqno{(4.2)}
$$
which coincides with $N$ if all $t_j$ are integers, 
 $t_j=\al_j, \, 1\le j\le m$ (cf. (2.1)).

\vspace{1ex} 

Further, introduce two functions 
$$ \tilde{\si} (t)= \tilde{\si}_{n,m} (t):= 
\prod_{j=1}^m \frac{p_j^{t_j+1}-1}{p_j-1}
 \ \cdot\prod_{j=n+1}^{m} \frac{p_j^3-1}{p_j-1}
\ \cdot\prod_{j=m+1}^{k} \frac{p_j^2-1}{p_j-1};
$$
$$ \tilde{G}(t) =  \tilde{G}_{n,m}(t) :=
 \frac{\tilde{\si}(t)}
{\tilde{N}(t)} \cdot \frac{1}{\log\log\tilde{N}(t)};
\qq \tilde{\eta}(t) = \tilde{\eta}_{n,m}(t):=\log\tilde{N}(t),
\eqno{(4.3)}
$$
which are the analogues of $\si(N), G(N), \eta$, cf. (1.5), (2.2).
Further, since 
$$ \frac{\tilde{\si}(t)}{\tilde{N}(t)}
=\prod_{j=1}^m \frac{p_j^{t_j+1}-1}{p_j^{t_j+1}} 
 \ \cdot
\prod_{j=m+1}^{n} \frac{p_j^3-1}{p_j^3}
\ \cdot
\prod_{j=n+1}^{k}  \frac{p_j^2-1}{p_j^2}
\cdot
 \prod_{j=1}^k \frac{p_j}{p_j-1};
\eqno{(4.4)}
$$
 (cf.  (4.2), (4.3)) one obtains  the following identity 
for the function
$$ \tilde{g}(t)=\tilde{g}_{n,m}(t):= \log \tilde{G}(t)  
= \sum_{j=1}^m \log\left(1-\frac{1}{p_j^{t_j+1}} \right)
- \log\log\tilde{\eta}(t)
 $$
$$   + \sum_{j=m+1}^n \log\left(1-\frac{1}{p_j^3}\right)
+ \sum_{j=n+1}^k \log\left(1-\frac{1}{p_j^2}\right)
 -  \sum_{j=1}^k \log\left(1-\frac{1}{p_j}\right) .
\eqno{(4.5)}
$$

By construction,   $\hat{g}(t)$ is bounded from above and tends to $0$ as $t\to\infty$; 

therefore there exists
 $\hat{g}_k = \hat{g}(n(N), m(N)) :=\max\{\tilde{g}(t): \, t\in ( \mathbb{R}_+)^m\}$.

Besides, it is obvious that  $\log G(N) \le \hat{g}_k$ for arbitrary 
 $N\in{\bf U_{1,k}}$ and 

 the same is valid also for all 
$N\in{\bf W_k}$ (since  $k\in{\bf E}$). 

\vspace{3ex} 

{\bf Lemma 2.} {\it For any $k\in\N, \, n, m$ as in (4.2),  the following upper estimate: 
 $$\hat{g}_k  < S_k - \log\log\theta_k
 - \frac{2\sqrt{2}-\delta_k}{\sqrt{p_k}\log p_k},\ 
\qq \delta_k \searrow 0.
\eqno{(4.6)}
$$

 holds true.}

\vspace{1ex} 

\noindent
$\triangleright$\ 
To find $\hat{g}_k$ we differentiate $\tilde{g}(t)$ with respect to   variables $t_j, j\in\{1,\dots , m\}$
$$ \frac{\partial\tilde{g}}{\partial t_j} 
=\left(\frac{1}{p_j^{t_j+1}-1}
 - \frac{1}{\tilde\eta \log\tilde\eta} \right) \log p_j.
\eqno{(4.7)}
$$

Equating  all these derivatives to $0$,
we come to the system of $(m+1)$

 equations with $(m+1)$ unknowns (cf. also (4.3)): 

\vspace*{\fill}

\clearpage

 $$ p_j^{t_j+1} = 1+\tilde{\eta} \log\tilde{\eta}, 
\ j\in\{1,\dots, m\};
$$
$$
 \q \tilde{\eta}=\sum_{j=1}^m t_j\log p_j + 2\sum_{i=m+1}^n \log p_i
+\sum_{i=n+1}^k \log p_i.
\eqno{(4.8)}
$$

From the first $m$ of these equations one obtains:
$$ p_j^{t_j}=\frac{1+\tilde{\eta}\log\tilde{\eta}}{p_j} 
\Rightarrow t_j\log p_j = \log (1+\tilde{\eta}\log\tilde{\eta}) - \log p_j;
\ 1\le j\le m.
\eqno{(4.9)}
$$

Summing here over $j$ from $1$ to $m$ one comes to:
\vspace{-1ex} 
$$ \sum_{j=1}^m t_j\log p_j 
= m\cdot\log (1+\tilde{\eta}\log\tilde{\eta}) - \theta_m. 
\eqno{(4.10)}
$$

\vspace{-1ex} 

But the left-hand side (according to the last equation (4.8) and (1.2)) 
equals  $\tilde{\eta}- 2(\theta_n-\theta_m) - (\theta_k -\theta_n)
 =\eta-\theta_k -\theta_n  +2\theta_m$
 and thus one comes
to the single equation for finding $\tilde{\eta}$:
$$ \tilde{\eta} 
= (\theta_k + \theta_n - 3\theta_m) + m\, \log (1+\tilde{\eta}\log\tilde{\eta}).
\eqno{(4.11)}
$$

Having found $\tilde{\eta}^*$ from this equation, one can calculate 
by means of (4.9) the (unique) stationary point 
$t^*=(t_1^*,\dots, t_m^*)$, 
 in which $\tilde{g}_k(t)$  attains its maximum,
 and then from (4.5) (cf. also the definition (1.3) of $S(x)$) one comes to:
$$ t_j^*= \frac{\log(1+\tilde{\eta}^*\log\tilde{\eta}^*)}{\log p_j} - 1;
$$
$$
 \hat{g}_k = S_k - \log\log\tilde{\eta}^* 
+ m \log\left(1-\frac{1}{1+ \tilde{\eta}^*\log\tilde{\eta}^*}\right) 
$$
$$+\sum_{i=m+1}^n \log\left(1-\frac{1}{p_i^3}\right)
+\sum_{i=n+1}^k \log\left(1-\frac{1}{p_i^2}\right).
\eqno{(4.12)}
$$

\noindent
{\bf 4.2.} In order to obtain good approximation for the root of the 
determining equation (4.11), we will need the following auxiliary assertion.

\vspace{1ex}

{\bf Lemma 3.} {\it For $A> 2B > 1$
 the equation 
\vspace{-1ex}
$$x=A+B\log (1+x\log x)
\eqno{(4.13)}
$$ 
has a unique root 
$x^*=x^*(A, B)>0$, which satisfies the estimates}:
\vspace{-1ex}
$$ A+B\log(1+A\log A) < x^* 
 < A+B\log(1+A\log A) \left(1+  \frac{2B}{A - 2B}\right)
\eqno{(4.14)}
$$ 

\vspace*{\fill}

\clearpage

\noindent
$\triangleright$ \ Consider a function $y(x):=A+B\log (1+x\log x)$. 
Under assumptions of Lemma 
$\min y(x) = y(1/e)=A+B\log(1-1/e)>B\log(e-1)>1$; therefore equation
(4.12) has no roots on $(0, 1]$. On the other hand, on the interval
 $(0,+\infty)$  the function $y(x)$ is increasing, {\it concave} and the ratio 
$y(x)/x \to 0,\\ x\to\infty;$ hence  there is exactly one root
 $x^*=x^*(A, B) $, which is a limit of iterations
$x_0=A, x_s:=y(x_{s-1})>x_{s-1}, \ s\in\N$. 

 The left-hand side is $x_1$. To obtain the upper estimate (3.14) we note that the mapping  $x\mapsto y(x)$ of $[A, \infty)$ into itself is a contraction with the coefficient
 $K<\sup\{ |y^{\,\prime}(x)|: x>A\}=B(1+\log A)/(1+A\log A)<2B/A<1$; therefore:
$$ x^*-x_1 < (x_1-x_0)\,\frac{K}{1-K}
<B\log(1+A\log A)\, \frac{2B}{A-2B};
 \eqno{(4.15)}
$$

\noindent
and thus the proof of Lemma 2  is complete $\Box.$

\vspace{1ex}
\noindent
{\bf 4.3.} Returning to the proof of Lemma 1,  
let us note that from definitions,  Theorem 2 (cf. (3.12))
 and Proposition 3 (cf. (3.6){\bf(ii)}) one immediately obtains
the relationships
$$ \theta_n\approx \sqrt{2\theta_k}, \q
 \theta_m \approx \sqrt[3]{3\theta_k},\ \ k\to\infty. 
 \eqno{(4.16)}
$$
 
On the other hand, according to Proposition 4 ((3.10){\bf(ii)})
 $m\approx \theta_m/\log\theta_m$.

Now apply Lemma 1 with (cf. (4.13))
$$ A:=\theta_k +\theta_n  -3\theta_m, 
\q B:=m\approx \theta_m/\log \theta_m
\approx 3\sqrt[3]{3\theta_k}/\log\theta_k, 
 \eqno{(4.17)}
$$

 whence 
 $B/A\approx \theta_m/ (\theta_k \log \theta_m)
\approx 3^{4/3}/(\theta_{\ k}^{2/3}\log \theta_k)$ follows,
and thus by virtue of (3.13), (3.15) one has for $k\to\infty$:
$$ \tilde{\eta}^*=\tilde{\eta}^*_{k,m,n}= (\theta_k+\theta_n -3 \theta_m)
+ 3^{4/3}\theta_{\ k}^{1/3} (1+o(1))
=\theta_k + \sqrt{2 \theta_k}  + O(\theta_{\ k}^{1/3}) . 
 \eqno{(4.18)}
$$

\vspace{1ex} 

Substituting these values of $m, n$ and $\tilde{\eta}^*$ 
into (4.11), one obtains

 (explanations below) 
\vspace{-1ex}
$$ \hat{g}_k=S_k - \log\log\theta_k 
- (\log\log\tilde{\eta}^* - \log\log\theta_k) 
+ \sum_{j=n+1}^k \log \left(1-\frac{1}{p_j^2}\right) 
$$
$$ +O\left(\frac{1}{\sqrt{\theta_k} \log^2\theta_k}\right)
 = S_k - \log\log\theta_k 
- \frac{2\sqrt{2}+O(1/\log p_k)}{\sqrt{p_k}\log p_k}.
 \eqno{(4.19)}
$$

\vspace{2ex} 

Here we have used the relationships:

\vspace*{\fill}

\clearpage

$$ \log\log\tilde{\eta}^* - \log\log\theta_k
\approx \frac{\theta_n}{\theta_k\log\theta_k} 
\approx \frac{\sqrt{2}}{\sqrt{p_k}\log p_k};
$$
$$ m \log\left(1-\frac{1}{1+ \tilde{\eta}^*\log\tilde{\eta}^*}\right)
\approx -\frac{3^{4/3}}{\theta_k^{2/3}\log^2\theta_k}
=O\left(\frac{1}{p_k^{2/3}\log^2 p_k}\right).
 \eqno{(4.20)}
$$
\vspace{1ex} 

and have also taken advantage of Proposition A,
 cf. (2.6) with $\la=2$.

\vspace{2ex} 

Hence it follows immediately, 
that two last summands  in (4.12) are:
$$ \sum_{i=m+1}^n \log\left(1-\frac{1}{p_i^3}\right)
=- \frac{1}{2p_m^{2}\log p_m}
+O\left(\frac{1}{p_m^{2}\log^2 p_m}\right)
=O\left(\frac{1}{p_k^{2/3}\log p_k}\right);
$$
$$
\sum_{i=n+1}^k \log\left(1-\frac{1}{p_i^2}\right)
= -\frac{1}{p_n\log p_n} 
+ O\left(\frac{1}{p_n\log^2 p_n} \right)
$$
$$=- \frac{\sqrt{2}}{\sqrt{p_k}\log p_k} 
+ O\left(\frac{1}{\sqrt{p_k}\log^2p_k} \right),
\eqno{(4.21)} 
 $$ 

\vspace{1ex} 

and this completes the proof of the asymptotic formula (4.19)
 for $\tilde{g}_k^*$, 

which in turn implies (4.6) $\Box$.
\vspace{1ex}

The proof of the part {\bf (II)} and of the whole Theorem 1  
is thus complete.

\vspace*{\fill}

\clearpage

\centerline{\bf LIST OF REFERENCES }

\vspace{3ex}

\noindent
[1] Mertens F. {\it $\ddot{U}ber$ einige asymptotische Gesetze der Zahlentheorie}.

 J. Reine Angew.  Math.,  {\bf 77}, 1874, pp. 289 -- 338.

\vspace{2ex}

\noindent
[2]\   Nicolas J.-L. \  {\it Petites valeurs de la fonction d'Euler.} 

 Journal of Number Theory, vol. 17, no. 3, 1983, p. 375 -- 388.

\vspace{2ex}

\noindent
[3]\  Gronwall T.H. {\it Some asymptotic expressions 
in the theory of numbers}. 

Trans. Amer.  Math. Soc. V. 14 (1913), pp. 113 -- 122.

\vspace{2ex}
\noindent
[4]\   Ramanujan S. \ {\it Highly composite numbers, annotated 
and with a foreword} 

{\it by J.-L. Nicolas and G. Robin},
\ Ramanujan J. V. 1,  1997, pp. 119 -- 152.

\vspace{2ex}

\noindent
[5]\   Robin G. \ {\it Grandes valeurs de la fonction somme 
des diviseurs et hypoth\`ese }

{\it de  Riemann.} J. Math. Pures Appl. V. 63 (1984), 
pp. 187 -- 212.

\vspace{2ex}

\noindent
[6]\    Caveney G., Nicolas J.-L., Sondow J. 
\ {\it Robin's theorem, primes,  }

{\it and a new elementary reformulation of Riemann Hypothesis.}

INTEGERS 11 (2011),  \#A33, pp 1 -- 10.

\vspace{2ex}

\noindent
[7]\    Caveney G., Nicolas J.-L., Sondow J.
 \ {\it On SA, CA, and GA numbers}

{\ https://arxiv.org/abs/1112.6010v2},  2012.

\vspace{2ex}

\noindent
[8]\ Kalyabin G.\,A. {\it Refinement of Mertens formula
and Robin inequality.}

 International Conference 
{\bf "P. Chebyshev Mathematical Ideas and 

Their Applications to Natural Sciences"},  commemorating 

the 200-th anniversary of P.L. Chebyshev, 
the great Russian mathematician. 

 Short Papers. DOI: 10.51790/chebconf-2021,
 p. 188 - 189.

\vspace{2ex}

\noindent
[9]\ Kalyabin G.\,A. {\it Remainder in modified Mertens formula

and Ramanujan inequality.}

{\ https://arxiv.org/abs/2201.02663},  2022.

\vspace{2ex}

\noindent
[10]\ Kalyabin G.\,A. {\it One-step $G$-unimprovable numbers}. 
 
{\ https://arxiv.org/abs/1810.12585},  2018.

\vspace{2ex}

\noindent
[11]\ Kalyabin G.\,A. {\it {\bf RH}-dependent 
estimates of remainder

\ in modified Mertens formula}.

{\ https://arxiv.org/abs/2205.05931},  2022.

\vspace{2ex}

\noindent
[12]\ Narkiewicz W. {\it The Development of Prime Number Theory}. 
 
\ \, {Springer-Verlag Berlin Heidelberg New York},  2000.

\end{document}